\documentclass[reqno]{amsart}
\usepackage{verbatim}
\usepackage{amssymb}
\usepackage{enumerate}

\newtheorem{theorem}{Theorem}[section]
\newtheorem{proposition}[theorem]{Proposition}
\newtheorem{lemma}[theorem]{Lemma}
\newtheorem{corollary}[theorem]{Corollary}
\newtheorem{claim}{Claim}[theorem]
\newtheorem{case}{Case}

\theoremstyle{definition}

\newtheorem{definition}[theorem]{Definition}

\theoremstyle{remark}
\newtheorem{remark}{Remark}

\newif\ifdeveloping

\ifdeveloping
\usepackage[notref,notcite]{showkeys}
\fi

\newcommand{\mc}[1]{\mathcal{#1}}

\newcommand{\mb}[1]{\mathbf{#1}}
\newcommand{\mf}[1]{\mathfrak{#1}}

\newcommand{\setm}{\setminus}
\newcommand{\empt}{\emptyset}
\newcommand{\subs}{\subset}

\newcommand{\oo}{{{\omega}_1}}
\newcommand{\dom}{\operatorname{dom}}

\def\<{\left\langle}
\def\>{\right\rangle}
\def\br#1;#2;{\bigl[ {#1} \bigr]^ {#2} }

\newcommand{\supp}{\operatorname{supp}}
\newcommand{\tip}{\operatorname{tp}}
\newcommand{\id}{\operatorname{id}}
\newcommand{\rank}{\operatorname{rank}}

\newcommand{\force}{\Vdash}

\def\lev#1;#2;{\operatorname{I}_{#1}(#2)}

\title{Wide scattered spaces and  morasses
}

\author[L. Soukup]{Lajos Soukup }
\address{Alfr{\'e}d R{\'e}nyi Institute of Mathematics, Hungarian
  Academy of Sciences}
\email{soukup@renyi.hu}
\thanks{The second author was partially supported by
Hungarian National Foundation for Scientific Research grants no.
61600 and 68262}

\date{\today}

\begin{document}

\begin{abstract}
We show that it is relatively consistent with ZFC that $2^\omega$ is arbitrarily large and
every sequence $\mb s=\<s_{\alpha}:{\alpha}<{\omega}_2\>$ of infinite
cardinals with $s_{\alpha}\le 2^\omega$
is the cardinal sequence of some
locally compact scattered space.
\end{abstract}

\maketitle

\section{Introduction}
If $X$ is a scattered topological space, and ${\alpha}$ is an ordinal,
denote by $\lev {\alpha};X;$ the ${\alpha}$th Cantor-Bendixson level
of $X$.
The {\em cardinal sequence} of  $X$, $\operatorname{SEQ}(X)$, is the
sequence of the cardinalities of the infinite  Cantor-Bendixson levels
of $X$, i.e.
\begin{displaymath}
\operatorname{SEQ}(X)=\bigl\langle\ |I_{\alpha}(X)|:{\alpha}<ht^-(X)\ \bigr\rangle,
\end{displaymath}
where $ht^-(X)$, the {\em reduced height} of $X$, is 
the  minimal ${\beta}$ such that  $\lev \beta;X;$ is finite.
If $\delta$ is an ordinal, we denote by ${\mathcal
C}(\delta)$ the class of all cardinal sequences of length $\delta$
of locally compact scattered (LCS, in short) spaces.

Let $\<{\kappa}\>_{\alpha}$  denote the
constant ${\kappa}$-valued sequence of length ${\alpha}$.

\begin{theorem}[Baumgartner, Shelah, \cite{BS}]
It is relatively consistent with ZFC that $\<{\omega}\>_{{\omega}_2}\in \mc C({\omega}_2)$.   
\end{theorem}
Refining their argument, first   Bagaria, \cite{Bag}, proved that 
${}^{{\omega}_2}\{{\omega},{\omega}_1\}\subs \mc C({\omega}_2)$ in some  ZFC model,
then Martinez and Soukup, \cite{MS-Universal}, showed that  
$2^{\omega}={\omega}_2$ and 
${}^{{\omega}_2}\{{\omega},{\omega}_1, {\omega}_2\}\subs \mc C({\omega}_2)$
is also consistent.

For a long time
$\omega_2$ was  a mystique barrier in both height and width. 
In this paper we can construct wider spaces.

\begin{theorem}\label{tm:main}
If GCH holds and ${\lambda}\ge {\omega}_2$ is a regular cardinal, then
in some cardinal preserving generic extension $2^{\omega}={\lambda}$
and 
every sequence $\mathbf s=\<s_{\alpha}:{\alpha}<{\omega}_2\>$ of infinite
cardinals with $s_{\alpha}\le {\lambda}$
is the cardinal sequence of some
locally compact scattered space.  
\end{theorem}

We will find the suitable generic extension in three steps:
\begin{enumerate}[(I)]
\item The  first extension adds a ``{\em strongly stationary 
strong $({\omega}_1,{\lambda})$-semimorass}'' to the ground model
(see Definition \ref{df:strong_semimorass} and Theorem \ref{tm:morass}). 
\item Using that strong semimorass the second extension adds
a {\em  $\Delta({\omega}_2\times\lambda)$-function} to the first
extension
(see Definition \ref{df:sDelta} and Theorem \ref{tm:sDelta}).
\item Using the  $\Delta({\omega}_2\times {\lambda})$-function
we add an ``{\em  LCS space with stem}'' to the second model and we show that 
those2 space alone guarantees that   every sequence 
$\mb s=\<s_{\alpha}:{\alpha}<{\omega}_2\>$ of infinite
cardinals with $s_{\alpha}\le {\lambda}$ is the cardinal sequence of some
locally compact scattered space (see Theorem \ref{tm:space}). 
\end{enumerate}
Steps  (I) and (II) are based on works of P. Koszmider,
see \cite{Ko1} and \cite{Ko_V}.

\section{Strong semimorasses}

If ${\rho}$ is a function  and $X$ is set, write
${\rho}[X]=\{{\rho}({\xi}):{\xi}\in X\}$.

If $X$ and $Y$ are sets of ordinals with $\tip(X)=\tip(Y)$,
denote the unique order preserving bijection between
$X$ and $Y$ by ${\rho}_{X,Y}$.

For $X\in \br {\lambda};{\omega};$ and  
$\mc F\subs \br {\lambda};{\omega};$
let $\mc F\restriction X=\{Y\in \mc F:Y\subsetneq X\}$.

If $X$, $X_1$ and $X_2$ are sets of ordinals, we write
\begin{enumerate}[]
\item  $X=X_1\oplus X_2$ iff $\tip(X_1)=\tip(X_2)$,
 $X=X_1\cup X_2$ and
${\rho}_{X_1,X_2}\restriction X_1\cap X_2= \id
$;
\item
$X=X_1\odot X_2$ iff $\tip(X_1)=\tip(X_2)$,
 $X=X_1\cup X_2$ and
$X_1\cap X_2<X_1\setm X_2<X_2\setm X_1$;
\end{enumerate}
and
\begin{enumerate}[]
\item $X=X_1\otimes X_2$ iff   $X=X_1\oplus X_2$
and $X\cap {\omega}_2=(X_1\cap {\omega}_2)\odot (X_2\cap {\omega}_2)$.
\end{enumerate}

In \cite{Ko1}  Koszmider introduced the notion of  semimorasses
and proved several properties 
concerning that structures. Unfortunately, in our proof  we need  structures 
with a bit stronger properties.

\begin{definition}\label{df:strong_semimorass}
Let ${\omega}_1\le {\lambda}$ be a cardinal.
A family $\mc F\subs \br {\lambda};{\omega};$ is a 
{\em strong $({\omega}_1,{\lambda})$-semimorass}
iff\begin{enumerate}[(M1)]
\item $\<\mc F,\subseteq\>$ is {\em well-founded},
(and so we have the {\em rank} function on $\mc F$),
\item $\mc F$ is {\em locally small}, i.e.
$|\mc F\restriction X|\le {\omega}$ for each $X\in \mc F$.
\item $\mc F$ is {\em homogeneous}, i.e. 
$\forall X, Y\in \mc F$ if $\rank(X)=\rank(Y)$
then $\tip(X)=\tip(Y)$ and
$\mc F\restriction Y=\{{\rho}_{X,Y}[Z]:Z\in \mc F \restriction X\}$. 
\item $\mc F$ is {\em directed}, i.e. $\forall X,Y\in \mc F$
($\exists Z\in \mc F$) $X\cup Y\subs Z$.
\item $\mc F$ is {\em strongly locally semidirected}, i.e. $\forall
  X\in \mc F$ either 
(a) or (b) holds: 
  \begin{enumerate}[(a)]
  \item $\mc F\restriction X$ is directed,
  \item $\exists X_1, X_2\in \mc F$ 
$\rank(X_1)=\rank(X_2)$,
$X=X_1\otimes X_2$, 
and $\mc F\restriction X=(\mc F\restriction X_1)\cup (\mc F\restriction
X_2)\cup \{X_1, X_2\}$.    
\end{enumerate}
\item $\mc F$ {\em covers} ${\lambda}$, i.e. $\cup \mc F={\lambda}$.
\end{enumerate}
If in (M5)(b)  we weaken the assumption
$X=X_1\otimes X_2$ to $X=X_1\oplus X_2$ 
then we obtain the definition
of an {\em  $({\omega}_1,{\lambda})$-semimorass} (see \cite[Definition 1]{Ko1}).
Moreover,
a strong $({\omega}_1,{\omega}_2)$-semimorass
is just Velleman's simplified $({\omega}_1,{\omega}_2)$-morass.
\end{definition}

\begin{definition}
A family $\mc F\subs \br {\lambda};{\omega};$
is {\em strongly stationary} iff 
for each function $c:\br \mc F;<{\omega};\to \br
{\lambda};{\omega};$
there are stationary many $X\in \mc F$ such that $X$ is {\em $c$-closed}, i.e. 
$c(X^*)\subs X$ for each $X^*\in \br \mc F\restriction X;<{\omega};$.  
\end{definition}

\begin{theorem}\label{tm:morass}
If  $2^{\omega}={\omega}_1< {\lambda}={\lambda}^{\oo}$ then
there is a ${\sigma}$-complete ${\omega}_2$-c.c. forcing notion 
$P$ such that 
$$\notag
V^{P}\models \text{ ``{${\lambda}^{\oo}={\lambda}$ and there is a strongly
    stationary 
strong $({\omega}_1,{\lambda})$-semimorass $\mc F$. 
}''}  
$$
\end{theorem}

We say that a family $p\subs \br {\lambda};{\omega};$
is {\em neat} iff
 $\cup p=\cup (p\setm \{\cup p\})$.

\begin{proof}[Proof of Theorem \ref{tm:morass}]

Define $P=\<P,\le\>$ as follows. Let 
\begin{equation}\notag
P=\{p\subs \br {\lambda};{\omega};:
|p|\le {\omega}, \cup p\in p, 
\text{$p$ is neat and satisfies (M1)--(M5)}\}.
\end{equation}
Write $\supp(p)=\cup p$ for $p\in P$.
Clearly  $\supp(p)$ is the  $\subs$-largest element of $p$.
Put
\begin{equation}\label{eq:lep}
p\le q\text{ iff } \supp(q)\in p\land q=(p\restriction \supp(q))\cup\{\supp(q)\}.  
\end{equation}
$P$ is  ${\sigma}$-complete. Indeed, if $p_0\ge p_1\ge p_2\dots $ then let
\begin{equation}\notag
p=\cup_{n<{\omega}}p_n \cup\{\cup_{n<{\omega}} \supp({p_n})\}. 
\end{equation}
Then $p\in P$ and $p\le p_n$ for each $n$.

\begin{definition}
We say that two conditions 
$p$ and $p'$ are {twins} iff
\begin{enumerate}[(i)]
\item $\tip (\supp (p))=\tip(\supp (p'))$,
\item $\supp(p)\cup \supp(p')=\supp(p)\otimes \supp(p') $,
\item $p'={\rho}_{\supp(p),\supp (p')}[p]$.
\end{enumerate}
\end{definition}

\begin{lemma}\label{cl:twin}
If $p$ and $p'$ are twins then 
they have a common extension in $P$ 
\end{lemma}

\begin{proof}
Write $D=\supp(p)$ and $D'=\supp(p')$.
Put $r=p\cup p'\cup\{D\cup D'\}$.
We show that $r$ is a common extension of $p$ and $p'$.

\noindent {\bf Claim:} {\em $p=(r\restriction D) \cup\{D\}$ and
  $p'=(r\restriction D')\cup\{D'\}$.}

Indeed, assume that $X\in r\restriction D$.
Then $X\in p$ or $X\in p'$.
If $X\in p'$ then $X\subs D'$, and so   $X\subs D\cap D'$.
Since ${\rho}_{D,D'}\restriction D\cap
D'=\id$ it follows that $X={\rho}^{-1}_{D,D'}[X]\in p$.
So $r\restriction D\subs p$, which proves the Claim.

\smallskip

First we check that 
 that 
 $r\in P$.
(M1) and (M2) are clear.
Since $\supp(r)=\supp(p)\cup \supp(p')$,
$r$ is neat. 
$r$ has the largest element 
$\supp(r)=D\cup D'\in r$,
and so (M4) also holds.
In (M5) we have just one  new instance $X=\supp(r)$.
But in this case the choice $X_1=D$ and $X_2=D'$ works.
 To check (M3) assume that  $X,Y\in r$, $\rank(X)=\rank(Y)$.
If $X,Y\in p$ or $X,Y\in p'$ then we can apply that 
$p$ and $p'$ satisfy (M3).
So we can assume that $X\in p\setm p'$ and 
$Y\in p'\setm p$.
Let $X'={\rho}_{D,D'}[X]$.  Then
$\rank(X')=\rank(X)=\rank(Y)$ and 
$X',Y\in p'$.
Since $p'$ satisfies (M3), we have 
 $\tip(X')=\tip(Y)$,  and so $\tip(X)=\tip(Y)$.
Since 
${\rho}_{X',Y}:p' \restriction X'\to p' \restriction Y$ 
is an isomorphism,  and ${\rho}_{X,Y}= {\rho}_{D,D'}\circ {\rho}_{X',Y} $
it follows that 
${\rho}_{X,Y}:p\restriction X\to  p'\restriction Y$ is also an 
isomorphism. However  $:p\restriction X=r\restriction X$ and 
$p'\restriction Y=r\restriction Y$ by the Claim,  and so 
${\rho}_{X,Y}:r\restriction X\to  r\restriction Y$ is also an 
isomorphism,
  which proves (M3).

Finally  $r\le p, p'$ follows immediately from the Claim.
\end{proof}

\begin{lemma}\label{cl:w2}
$P$ satisfies ${\omega}_2$-c.c.  
\end{lemma}

\begin{proof}
Assume that $\{r_{\alpha}:{\alpha}<{\omega}_2\}\subs P$.
Write $D_{\alpha}=\supp({r_{\alpha}})$ for ${\alpha}<{\omega}_2$.
By standard argument we can find $I\in \br \omega_2;\omega_2;$
such that

\begin{enumerate}[(a)]
\item $\{D_{\alpha}:{\alpha}\in I\}$ forms a $\Delta$-system
with kernel $D$, and $\tip(D_\alpha)=\tip(D_\beta)$ for $\alpha,\beta\in I$,		
\item For ${\alpha}<{\beta}\in I$ we have 
$D\cap {\omega}_2<D_{\alpha}\setm D<D_{\beta}\setm D$,
\item ${\rho}_{D_{\alpha},D_{\beta}}[D_{\alpha}\cap
  {\omega}_2]=D_{\beta}\cap {\omega}_2$ and 
${\rho}_{D_{\alpha},D_{\beta}}\restriction D=\id$,
\item $r_{\beta}=\{{\rho}_{D_{\alpha},D_{\beta}}[X]:X\in r_{\alpha}\}$.
\end{enumerate}
Then for each 
 ${\alpha}\ne {\beta}\in I$  the conditions
$r_{\alpha}$ and $r_{\beta}$ are twins, so they are compatible by 
Lemma \ref{cl:twin}.
\end{proof}

\begin{lemma}
(a)  $\forall {\alpha}\in {\omega}_2$
  \begin{equation}\notag
  D_{\alpha}=\{p\in P: \supp(p)\cap ({\omega}_2\setm
  {\alpha})
\ne \empt\}  
  \end{equation}
is dense in $P$.\\
(b) $\forall {\beta}\in {\lambda}\setm {\omega}_2$
\begin{equation}\notag
E_{\beta}=\{p\in P:{\beta}\in \supp(p)\}  
\end{equation}
is dense in $P$.
\end{lemma}

\begin{proof}
(a) For each $q\in P$  and $ {\alpha}<{\omega}_2$
there is
 $q'$  such that  $q$ and $q'$  are twins and 
$\supp(q')\cap ({\omega}_2\setm {\alpha})\ne \empt$.
Then $q$ and $q'$ has a common extension $p\in D_{\alpha}$
by Lemma \ref{cl:twin}.\\
(b)   For all $q$ and ${\beta}\in {\lambda}\setm {\omega}_2$
there is $q'$ such that  $q$ and $q'$  are twins and 
${\beta}\in \supp(q')$. Then the common extension $p$ 
of $q$ and $q'$ is in $E_{\beta}$.
\end{proof}

Let $\mc G$ be a $P$-generic filter over $V$, and put
 $\mc F'=\cup\mc{G}$ and $F=\cup\mc F'$.
Then $\mc F'\subs \br \lambda;\omega;$ and so $F\subs \lambda$.
By the previous  lemma,
$F\supset  {\lambda}\setm {\omega}_2$
and $|F\cap {\omega}_2|={\omega}_2$.
So $\mc F=\{\rho_{F,\lambda}[X]:X\in \mc F'\}$ is a
strong $({\omega}_1,{\lambda})$-semimorass. 
To complete the proof of \ref{tm:morass} it is enough to prove the
following lemma.

\begin{lemma}\label{lm:ss}
$V^P\models$ ``$\mc F$ is strongly stationary''.  
\end{lemma}

\begin{proof}
It is enough to prove that if 
\begin{equation}\notag
q\force \dot {\mc C}\subs \br F;{\omega};
\text{ is club and }
\dot c:\br \mc F';<{\omega};\to \br F;{\omega}; 
\end{equation}
then there are  $p\le q$  and $C\in \br \lambda;\omega;$ such that 
$p\force $  ``$\check C\in \mc F'\cap \dot{\mc C} $ is $\dot c$-closed''. 

First we need a claim. 
\begin{claim}
If $p\force \check A\in \br F;{\omega};$
then 
$\exists p'\le p$ such that $A\subs \supp({p'})$.  
\end{claim}

\begin{proof}
If ${\alpha}\in A$ and 
$p\force \check{\alpha}\in F$ then there are $p'\le p$ and $p''\in P$
such that ${\alpha}\in \supp(p'')$ and $p'\force p''\in F$.
Then $p'$ and $p''$ have a common extension $q$, and 
then  ${\alpha}\in \supp(q)$ and $q\force  \check {\alpha}\in F$.

Since $P$ is ${\sigma}$-complete and $A$ is countable,  we are done by a 
straightforward induction.
\end{proof}
We will
choose a decreasing sequence 
$\<p_n:n<{\omega}\>\subs P$ and an increasing sequence
$\<C_n:n<{\omega}\>\subs \br {\lambda};{\omega};$
as follows.
Let $C_0=\empt$ and $p_0=q$.
If $p_n$ and $C_n$ are given, let 
$Z_n\supset {C_n}\cup \supp(p_n)$  and $p_n'\le p_n$ s.t. 
\begin{equation}\notag
p_n'\force \cup\{\dot c(X):X\in \br \check p_n;<{\omega};\}\subs Z_n\in \br
F;{\omega};.  
\end{equation}
Let  $p_{n+1}\le p'_n$ and $C_{n+1}\supset Z_n\cup C_n$ such that 
$C_{n+1}\subs \supp(p_{n+1})$ and
$p_{n+1}\force \check{C}_{n+1}\in \dot{\mc C}$.

Having constructed the sequence  finally put
$C=\cup\{C_n:n<{\omega}\}$ and
$p=\cup_{n<{\omega}} p_n\cup\{C\}$. 
Then $p\in P$, $p\le q$, $C=\supp (p)$.
Since $p\force $ ``{\em $\check C_n\in \dot {\mc C}$ 
and $\dot {\mc C}$ is club}'', we have $p\force $ ``$\check C\in \dot {\mc C}$.
Since  
$p\force $ $\dot c''\bigl[\br p_n;<\omega;\bigr]\subs \supp(p_{n+1})$,
we have  $p\force $ {\em $\dot C$ is $\dot c$-closed.}
Moreover $p\force p\subs \mc F'$, so 
$p\force \check C\in \mc F'$.

Putting these together we obtain that 
$p$ and $C$ have the desired properties, 
which proves the lemma.
\end{proof}
Since ${\left({\lambda}^{\oo}\right)}^{V[\mc F]}\le 
{\left((|P|+{\lambda})^{\oo}\right)}^{V}=
{\left({\lambda}^{\oo}\right)}^{V}=
{\lambda}$,
the proof of Theorem \ref{tm:morass} is complete.
\end{proof}

Next we investigate some properties of strong semimorasses.

\begin{lemma}\label{lm:main_lemma_7}
Let $\mc F\subs \br {\lambda};{\omega};$  be a  strongly stationary
strong $({\omega}_1,{\lambda})$-semimorass.
\begin{enumerate}[(1)]
\item If 
$X,Y\in \mc F$, $\rank(X)=\rank(Y)$, ${\alpha}\in X\cap Y\cap {\omega}_2$,
then 
$X\cap {\alpha}=Y\cap {\alpha}$.
\item If $X,Y\in \mc F$,  $X\subs Y$ and  $\rank(X)<{\alpha}<\rank(Y)$
then  there is $Z\in \mc F$ such that  
$\rank(Z)={\alpha}$ and $X\subs Z\subs Y$.
\item If $X\in \mc F$ and  $\rank(X)<{\alpha}<{\omega}_1$
then  there is $Z\in \mc F$ such that  
$\rank(Z)={\alpha}$ and $X\subs Z$.
\item If 
$X,Y\in \mc F$, $\rank(X)\le \rank(Y)$, 
and ${\alpha}\in X\cap Y\cap {\omega}_2$,
then 
$X\cap {\alpha}\subs Y\cap {\alpha}$.
\end{enumerate}
\end{lemma}

\begin{proof}
(1)
We prove the statement by induction on the minimal 
rank of    
$Z\in \mc F$ with $Z\supset X\cup Y$.

If  rank of $Z$ is minimal, then clearly
 $Z=Z_1\otimes Z_2$ where  $X\subs Z_1$ and $Y\subs Z_2$.
Let $X'={\rho}_{Z_1,Z_2}[X]\in \mc F\restriction Z_2$.
Since ${\alpha}\in Z_1\cap Z_2\cap {\omega}_2$,
 we have $Z_1\cap {\alpha}=Z_2\cap {\alpha}$ and so 
${\rho}_{Z_1,Z_2}\restriction ({\alpha}+1)=\id$.
Thus $X'\cap {\alpha}=X\cap {\alpha}$ and ${\alpha}\in X'.$
Since 
$X',Y\in \mc F\restriction Z_2$, ${\alpha}\in X'\cap Y$ and
$\rank(Z_2)<\rank(Z)$, by the inductive assumption we have
$X'\cap {\alpha}=Y\cap {\alpha}$.
Thus 
$X\cap {\alpha}=Y\cap {\alpha}$.

\smallskip

\noindent (2)
Easy by straightforward induction on $\rank(Y)$.

\smallskip

\noindent (3)
By straightforward induction on ${\alpha}$ there is $Y\in \mc F$
such that $X\subs Y$ and $rank(Y)\ge {\alpha}$. Then apply (2).  

\smallskip

\noindent (4)
By (3)
there is $Y'\supset X$ such that $\rank(Y)=\rank(Y')$.
Then apply (1) for $Y$ and $Y'$.
\end{proof}

In \cite{Ko_V} Koszmider proved several statements
for Velleman's simplified morasses. Here we need similar results 
for strong semimorasses. 
The following lemma corresponds to 
\cite[Fact 2.6-Fact 2.7]{Ko_V}. 

\begin{lemma}\label{lm:fact2.6}\label{lm:fact2.7}
\label{lm:fact2.8}\label{lm:fact2.9}
Let $\mc F\subs \br {\lambda};{\omega};$  be a  strongly stationary
strong $({\omega}_1,{\lambda})$-semimorass.
Assume $\lambda^\omega=\lambda$, fix an injective function  
$c:\mc F\to {\lambda}$, and consider the stationary set 
\begin{equation}\label{eq:fprime}
\mc F'=\{X\in \mc F:c(X^*)\in 
X\text{ for each $X^*\in \mc F\restriction X$}\}.  
\end{equation}
Assume that  
 $\mc F,\mc F',c\in M\prec H(\theta)$, 
$|M|={\omega}$, and $M\cap {\lambda}\in \mc F'$.
Then 
\begin{enumerate}[(1)]
\item 
 $\mc F\restriction M\cap {\lambda}\subs M$.
\item  $\rank(M\cap {\lambda})=M\cap {\omega}_1$.
\item
 If $Y\in \mc F$ with 
$\rank(Y)<{\delta}=M\cap {\omega}_1$
then there is $Z\in M\cap \mc F$ such that 
 $(M\cap {\lambda})\cap Y\subs Z$,
and  $\rank(Z)=\rank(Y)$.
\item If $\mc A\in \br \mc F;<{\omega};$
then there is $Z\in \mc F\cap M$
such that 
\begin{equation}\label{eq:Z}
\cup\{X\cap M:X\in \mc A, \rank(X)<M\cap \oo\}\subs Z.  
\end{equation}
\end{enumerate}
\end{lemma}

\begin{proof}
(1) Let $X\in \mc F\restriction M\cap {\lambda}$,
i.e. $X\in \mc F$ and $X\subsetneq M\cap {\lambda}$.
Then  
$X\subsetneq M\cap {\lambda}\in \mc F'$ implies
${\alpha}=c(X)\in M\cap {\lambda}$.
But $c,{\alpha}\in M$ and $c$ is injective, 
so $X=c^{-1}\{{\alpha}\}\in M$.

\noindent
(2)  If $X\subsetneq M\cap {\lambda}$, $X\in \mc F$,
then $X\in M$ by (1)  
 and so
 $\rank(X)\in M\cap {\omega}_1$. Thus
$\rank(M\cap {\lambda})\le M\cap {\omega}_1$.

Assume that ${\alpha}<M\cap {\omega}_1$.
Then
\begin{equation}\notag
M\models ``\text{$\exists X\in \mc F$ $\rank(X)={\alpha}$.}''  
\end{equation}
Thus there is $X\in M\cap \mc F$ such that $\rank(X)={\alpha}$.
Hence $\rank(M\cap {\lambda})\ge M\cap {\omega}_1$.

\noindent (3)
There is $Y'\supset Y$, ${Y'}\in \mc F$ and 
$\rank({Y'})=\rank(M\cap \lambda)$.
Let $Z={\rho}_{{Y'},M\cap \lambda}[Y]$.
Since $Y\cap (M\cap \lambda)\subs  Y'\cap (M\cap \lambda) $
and ${\rho}_{Y',M\cap \lambda}\restriction Y'\cap (M\cap \lambda) =\id$,
we have $Z\supset Y\cap (M\cap \lambda)$.  

\noindent (4)
Just apply (3) and the fact  that 
$\mc F$ is directed.  
\end{proof}

\section{A $\Delta({\omega}_2\times {\lambda})$-function}
Let $\lambda\ge \omega_2$ be an infinite cardinal and let $\pi:{\omega}_2\times {\lambda}\to {\omega}_2$ be the natural
projection: 
$\pi(\<{\xi},{\alpha}\>)={\xi}$.

\begin{definition}\label{df:sDelta}
(1) Assume that  $f$ is a function such that   
$\dom(f)\subs \br {\omega}_2\times {\lambda};2;$ and 
$f\{x,y\}\in \br {\pi}(x)\cap {\pi}(y);<{\omega};$ for each 
$\{x,y\}\in \dom(f)$.
We say that two finite subsets $d$ and $d'$ of ${{\omega}_2}\times
 {\lambda}$ are  {\em good for $f$} provided 
$\br d\cup d';2;\subs \dom(f)$ and
$\forall x\in d'\setm d$
$\forall y\in d\setm d'$
$\forall z\in d\cap d'\cap ({\omega}_2\times
{\omega})$
\begin{enumerate}[(S1)]
\item if ${\pi}(z)<{\pi}(x),{\pi}(y)$
then ${\pi}(z)\in f\{x,y\}$,
\item if ${\pi}(z)<{\pi}(y)$
then $f\{x,z\}\subs f\{x,y\}$,
\item if ${\pi}(z)<{\pi}(x)$
then $f\{y,z\}\subs f\{x,y\}$.
\end{enumerate}
(2)
 A function $f:\br {\omega}_2\times {\lambda};2;\to \br
 {\omega}_2;<{\omega};$
is a {\em $\Delta({\omega}_2\times {\lambda})$-function}
iff $f\{x,y\}\subs \min({\pi}(x),{\pi}(y))$
and for each sequence 
$\{d_{\alpha}:{\alpha}<{\omega}_1\}\subs \br
{\omega}_2\times{\lambda};<{\omega};$ there are
${\alpha}\ne {\beta}$ such that 
$d_{\alpha}$ and $d_{\beta}$ are good for $f$.
\end{definition}

\begin{remark}
The assumption $|f\{x,y\}|<{\omega}$,
instead of the usual $|f\{x,y\}|\le {\omega}$,
is not a misprint.  
\end{remark}

\begin{theorem}\label{tm:sDelta}
If $2^{\omega}=\oo<{\lambda}={\lambda}^{\oo}$ and there is a 
strongly stationary  
strong $({\omega}_1,{\lambda})$-semimorass, 
then  in some cardinal preserving generic extension
${\lambda}^{\oo}={\lambda}$
and 
 there is a 
$\Delta({\omega}_2\times {\lambda})$-function.
\end{theorem}

\begin{proof}
To start with fix a strongly stationary 
strong $({\omega}_1,{\lambda})$-semimorass
$\mc F\subs \br {\lambda};{\omega};$.
We can assume that 
\begin{equation}
\label{omega} \text{${\omega}\subs X$ for each $X\in \mc F$.}  
\end{equation}
Fix   an injective function  
$c:\mc F\to {\lambda}$, and consider the stationary set
\begin{equation}\label{eq:fprime2}
\mc F'=\{X\in \mc F:c(X^*)\in X\text{ for each $X^*\in \mc F\restriction X$}\}.  
\end{equation}

\begin{definition}\label{df:function_poset}
We define a poset $P=\<P,\le\>$ as follows:
$P$ consists of triples $p=\<a,f,\mc A\>$,
where $a\in \br {\omega}_2\times {\lambda};<{\omega};$,
$f:\br a;2;\to \mc P ({\pi}[a])$
with $f\{s,t\}\subs \min({\pi}(s), {\pi}(t))$,
$\mc A\in \br \mc F;<{\omega};$
such that 
\begin{equation}\label{eq:pdef}
\text{$\forall s,t\in a$ $\forall X\in \mc A$
if $s,t\in X\times X$ then 
$f\{s,t\}\subs X$.
}
\end{equation}
Write $p=\<a_p, f_p,\mc A_p\>$ for $p\in P$.
Put $p\le q$ iff $a_p\supset a_q$, $f_p\supset f_q$
and $\mc A_p\supset \mc A_q$.

For $p\in P$
let
\begin{equation}\notag
\supp(p_\nu)=\cup a_\nu\cup \cup \mc A_\nu.
\end{equation}
Clearly $\supp(p)\in \br {\lambda};\le{\omega};$.

If $\rho$ is a function and $x=\<a,b\>\in \dom(f)^2$, let
$\bar \rho(x)=\<\rho(a), \rho(b)\>$.
We say $p,q\in P$ are {\em twins} iff 
\begin{enumerate}[(A)]
 \item $\supp(p)$ and $\supp(q)$ have the same order type,
\item the unique $<_{On}$-preserving bijection $\rho$ between 
  $\supp(p)$ and $\supp(q)$ gives an {\em isomorphism} between $p$  and $q$, i.e.
\begin{enumerate}[(a)]
 \item $a_q=\bar \rho[a_p]$,
\item $\{\rho[X]:X\in \mc A_p\}=\mc A_q$,
\item for each $\{s,t\}\in \br a_p;2;$, $\rho[f_p\{s,t\}]=f_q\{\bar \rho(s), \bar \rho(t)\}$.
\end{enumerate}
 \end{enumerate}
\end{definition}

\begin{definition}[{\cite[Definition 22]{Ko1}}]
Let $\mc K\subs \br {\lambda};{\omega};$.
A poset $P$ is {\em $\mc K$-proper}
iff for some large enough regular cardinal $\theta $
if 
$M$ is a countable elementary submodel of 
$\mc H(\theta)$ with $P\in M$ and $M\cap {\lambda}\in \mc K$
then for each $p\in M\cap P$ there is an $(M,P)$-generic $q\le p$.  
\end{definition}

\begin{lemma}[{\cite[Fact 23]{Ko1}}]
If $\mc K\subs \br {\lambda};{\omega};$ is stationary and 
a poset $P$ is {$\mc K$-proper}, then forcing with $P$
preserves $\oo$.  
\end{lemma}

\begin{definition}[{\cite[Definition 24]{Ko1}}]
Assume that $P$ is a poset, $M\prec \mc H(\theta)$, $|M|={\omega}$,
$q\in P$, and $P, {\pi}_1,\dots, {\pi}_n\in M$.
We say that {\em the formula $\Phi(x,{\pi}_1,\dots, {\pi}_n)$ well-reflects
  $q$ in $M$} 
iff   
\begin{enumerate}[(1)]
\item $\mc H(\theta)\models \Phi(q, {\pi}_1,\dots, {\pi}_n)$,
\item  if $s\in M\cap P$ and $M\models \Phi(s,{\pi}_1,\dots {\pi}_n)$
then $q$ and $s$ are compatible in $P$.  
\end{enumerate}
\end{definition}

\begin{definition}[{\cite[Definition 25]{Ko1}}]
Assume that $P$ is a poset, $\mc K\subs \br {\lambda};{\omega};$.
We say that $P$ is {\em simply $\mc K$-proper} if the following holds:
for some/each large enough regular cardinal 
$\theta$ \\
IF
\begin{enumerate}[(i)]
\item $M\prec \mc H(\theta)$, $|M|={\omega}$,
\item $p\in P$, $P,p,\mc K\in M$,
\item $M\cap {\lambda}\in \mc K$,  
\end{enumerate}
THEN there is $p_0\le p$ such that for each  $q\le p_0$
  some formula $\Phi(x,{\pi}_1,\dots,  {\pi}_n)$ well-reflects
  $q$ in $M$.
\end{definition}

By lemmas \cite[Fact 23 and Lemma 26]{Ko1}
we have
\begin{lemma}
If $\mc K\subs \br {\lambda};{\omega};$ is stationary and 
a poset $P$ is { simply $\mc K$-proper}, then forcing with $P$
preserves $\oo$.  
\end{lemma}

To show that  $\omega_1$ is preserved we prove the following lemma.
\begin{lemma}\label{lm:simply}
  $P$ is simply $\mc F'$-proper.  
  \end{lemma}

Actually we will prove some stronger statement.
To formulate it we need some preparation.

If  $M\prec \mc H(\theta)$, $p\in P\cap M$, $M\cap \lambda\in \mc F$ and $q\in P$ let 
\begin{equation}\notag
p^{M}= \<a_p, f_p, \mc A_p\cup\{M\cap \lambda\}\>
\end{equation} 
and
\begin{equation}\notag
q\restriction M=\<a_q\cap M, f_q\restriction M, \mc A_q\cap M\>.
\end{equation} 

\begin{lemma}\label{lm:rest}
(1) If  $M\prec \mc H(\theta)$ and  $p\in P\cap M$ 
then $p^M\in P$. (2) If $q\le p^M$ then $q\restriction M\in P\cap M$
as well.
\end{lemma}

\begin{proof}
(1)  We should only check (\ref{eq:pdef}) for $p^M$.
Assume that $s,t\in a_p$ and $X\in  \mc A_p\cup\{M\cap \lambda\}$.
Since $p\in P$, we can assume $X=M\cap \lambda$.
However $s,t\in M$, and so $f_p\{s,t\}\in M$ as well by $p\in M$.
Since $|f_p\{s,t\}|\le \omega$, it follows $f_p\{s,t\}\subs M\cap \lambda=X$
which was to be proved.  \\
(2) It is straightforward that $q\restriction M\in P$.
To show $q\restriction M\in M$ we should check that
$f_q\restriction M\in M$. So assume that $s,t\in a_q\cap M$.
Then $s,t\in (M\cap \lambda)\times (M\cap \lambda)$ and 
$M\cap\lambda\in \mc A_{p^M}\subs \mc A_q$. So, by (\ref{eq:pdef}),
$f_q\{s,t\}\subs M\cap \lambda$. Since $f_q\{s,t\}$ is finite,
we have $f_q\{s,t\}\in M$.
\end{proof}

\begin{lemma}\label{cl:Phi}
 Assume that  $M\prec \mc H(\theta)$, $|M|={\omega}$, $P, \mc F\in M$, 
 $p\in P\cap M$,  
 $M\cap {\lambda}\in \mc F'$. 
Let ${\delta}=\rank(M\cap {\lambda})=M\cap {\omega}_1$.
Assume that
 $Z\in M\cap \mc F$ such that 
\begin{equation}\label{eq:fi}
Z\supset \cup\{X\cap M:X\in \mc A_q, \rank(X)<{\delta}\}.  
\end{equation}
Let $\Phi(x,Z, q\restriction M)$ be the following formula:
\begin{equation}\label{eq:Phi}
\text{`` $x\in P$, 
 $x\le q\restriction M$, $(a_x\setm a_{q\restriction M})\cap
(Z\times Z)=\empt.$''  }  
\end{equation}
Then \\
(1)
$\Phi(q,Z,q\restriction M)$ holds. \\
(2)
If $s\in M$,  
$M\models \Phi(s,Z,q\restriction M)$
and 
\begin{equation}\notag
h:[a_s\setm a_{q\restriction M}, a_q\setm a_{q\restriction
    M}]\to
\mc P( {\pi}[a_s\cup a_q])  
\end{equation}
such that 
\begin{equation}\label{eq:h}
h\{x,y\}\subs \min({\pi}(x), {\pi}(y))\cap
\bigcap\{X\in A_q:x,y\in X\times X, \rank(X)\ge {\delta}\},  
\end{equation}
then 
$r=\<a_s\cup a_q,f_s\cup f_q\cup h, \mc A_s\cup \mc A_q\>\in P$
is a common extension of $q$ and $s$.\\
(3)
$\Phi(x,Z,q\restriction M)$ well reflects $q$ in $M$. 
\end{lemma}

\begin{proof}
(1) Since $q\restriction M\in P$ by Lemma \ref{lm:rest}(2), 
we have $q\le q\restriction M$ by the definition of the relation $\le$.
Since $Z\in M$, we have $Z\times Z\subs M\times M\subs M$ and 
$a_q\setm  a_{q\restriction M}=a_q\setm M$.\\
(2)  
To show that $r\in P$ 
we need to check that $r=\<a_r, f_r, \mc A_r\>$ satisfies (\ref{eq:pdef}).
So let $x,y\in a_r$, $X\in \mc A_r$.
\begin{case}
$x,y\in a_s$ and $X\in \mc A_s$ or
$x,y\in a_q$ and $X\in \mc A_q$.  
\end{case}
Then everything is fine, because $s,q\in P$.
\begin{case}
$\{x,y\}\in \br a_s;2;$, $x\in a_s\setm a_q$ and $X\in \mc A_q$.  
\end{case}
If $\rank(X)<{\delta}$ then 
$(a_s\setm a_q)\cap (X \times X)=\empt$ by (\ref{eq:fi}), so 
$x\notin X\times X$. Thus (\ref{eq:pdef}) is void.

If $\rank(X)\ge {\delta}$ then ${\nu}=\min({\pi}(x), \pi(y))\in M\cap X$,
so $M\cap {\nu}\subs X\cap {\nu}$ by Lemma \ref{lm:main_lemma_7}(4).
Thus $f_r\{x,y\}=f_s\{x,y\}\cap {\nu}\subs M\cap {\nu}\subs X$.
\begin{case}
 $\{x,y\}\in \br a_q;2;$, $x\in a_q\setm a_s$ and $X\in \mc A_s$.
\end{case}
Since $(a_q\setm a_s)\cap M=\empt$, it is not possible that $X\in \mc A_s$.
Then 
$x\in a_q\setm a_s=a_q\setm a_{q\restriction M}$, so $x\notin M$.
However $X\subs M$ and so $x\notin X\times X$, so 
(\ref{eq:pdef}) is void.
\begin{case}
$x\in a_q\setm a_s $ and $y\in a_s\setm a_q$.  
\end{case}
Then
the assumption concerning $h$ in (\ref{eq:h})
 is stronger than (\ref{eq:pdef}).
Indeed, if  $y\in a_s\setm a_q$ then 
$y\notin Z\times Z$. So if  $y\in X\times X$ for some $X\in \mc A_q$
then  $\rank(X)\ge {\delta}$.

\medskip
\noindent
(3)
Define the function
\begin{equation}\notag
h:[a_s\setm a_{q\restriction M}, a_q\setm a_{q\restriction
    M}]\to
\br {\omega}_2;<{\omega};  
\end{equation}
by $h\{x,y\}=\empt$.
Then (\ref{eq:h}) holds, so $s$ and $q$ are comparable by (2),
which 
was to be proved.
\end{proof}

\begin{proof}[Proof of Lemma \ref{lm:simply}]
We can apply Lemma \ref{cl:Phi}  because by Lemma \ref{lm:fact2.9} we can pick
  $Z\in M\cap \mc F$ such that
 $Z\supset \cup\{X\cap M:X\in \mc A_q, \rank(X)<{\delta}\}.$  
\end{proof}



\begin{lemma}\label{lm:w2-cc}
$P$ satisfies ${\omega}_2$-c.c.  
\end{lemma}

\begin{proof}
Let $\{p_{\nu}:{\nu}<{\omega}_2\}\subs P$.
Put $p_\nu= \<a_\nu, f_\nu, \mc A_\nu\>$. Recall that 
$\supp(p_\nu)=\cup a_\nu\cup \cup \mc A_\nu$.
We can assume that 
\begin{enumerate}[(i)]
\item $\{\supp(p_{\nu}):{\nu}<{\omega}_2)\}$ forms 
a $\Delta$-system with kernel $D$
\item the conditions are pairwise twins
witnessed by functions ${\rho}_{{\nu},{\mu}}:\supp(p_{\nu})\to
\supp(p_{\mu})$.
\end{enumerate}
Fix $\nu<\mu<\omega_2$. Define the function $e$ as follows: 
\begin{equation}\notag
\dom(e)=[a_{\nu}\setm a_{\mu}, a_{\mu}\setm  a_{\nu}]
\text{ and $e\{s,t\}=\empt$.}
\end{equation}
We claim that 
\begin{equation}\label{eq:r}
r=\<a_{\nu}\cup a_{\mu}, f_{\nu}\cup f_{\mu}\cup e,
\mc A_{\nu}\cup \mc A_{\mu}\>  
\end{equation}
is a common extension of $p_{\nu}$ and $p_{\mu}$. We need to show that
$r\in P$. Since $p_\nu$ and $p_\mu$ are twins, we should check only
(\ref{eq:pdef}). 
So let  $t,s\in a_{\nu}\cup a_\mu$ and $X\in \mc A_\nu \cup \mc A_{\mu}$ with  
$s,t\in X\times X$. We can assume e.g. $X\in \mc A_\nu$.
Since $X\subs \supp(p_\nu)$, 
we have $s,t \in \supp(p_\nu)\times \supp(p_\nu)$.
Then  $ \supp(p_\nu)\times \supp(p_\nu)\cap a_\mu\subs a_\nu$ because
$ \supp(p_\nu)\times \supp(p_\nu)\cap a_\mu\subs D$ and $a_\nu, a_\mu$ are twins. 
Thus  we have
$s,t\in a_\nu$, and so   we are done because $p_\nu$
 satisfies   (\ref{eq:pdef}).
%
%
\end{proof}

To complete the proof of 
Theorem \ref{tm:sDelta}
we claim that if $\mc G$ is a $P$-generic filter then the function
\begin{equation}\label{eq:fdef}
f=\cup\{f_p:p\in \mc G\}  
\end{equation}
is a  $\Delta({\omega}_2\times {\lambda})$-function.

Assume that 
\begin{equation}\notag
p\force \{\dot d_{\xi}:{\xi}<{\omega}_1\}\subs 
\br {\omega}_2\times {\lambda};<{\omega};.  
\end{equation}
We can assume
\begin{equation}\notag
p\force \{\dot d_{\xi}:{\xi}<{\omega}_1\}
\text{ is a $\Delta$-system with kernel $\check d$}.
\end{equation}
Assume  $M\prec H(\theta)$, $|M|=\omega$,  
$p, \mc F'\<\dot d_{\xi}:{\xi}<{\omega}_1\>\in M$ and 
$X_0=M\cap {\lambda}\in \mc F'$.
Let
\begin{equation}\notag
p^M=\<a_p, f_p, \mc A_p\cup\{X_0\}\>.  
\end{equation}
Let $q\le p^M$, ${\xi}_1<{\omega}_1$ and 
$e_1\in \br {\omega}_2\times {\lambda};<{\omega};$  that 
\begin{equation}\notag
q\force \dot d_{{\xi}_1}=\check e_1 
\land
(e_1\setm d) \cap M=\empt \land  e_1\subs a_q.  
\end{equation}
Put $\delta=\rank(M\cap\lambda)$.
By Lemma \ref{lm:fact2.9} we can pick
 $Z\in M\cap \mc F$ such that 
\begin{equation}\notag
Z\supset \cup\{X\cap M:X\in \mc A_q, \rank(X)<{\delta}\}.  
\end{equation}
Consider the following formula
$\Psi(x,\xi,e)$ with free variables $x$, $\xi$ and $e$ and parameters
$Z,q\restriction M, \<\dot d_\xi:\xi<\oo\>, d \in M$:
\begin{equation}\notag
\Phi(x,Z,q\restriction M)\land
{\xi}\in {\omega}_1\land
 (x\force \dot d_{\xi}=e) \land e\subs a_x \land (e\setm d)\subs a_x\setm a_{q\restriction M} 
\end{equation}
where the formula $\Phi(x,Z,q\restriction M)$
was defined in (\ref{eq:Phi}) in Lemma \ref{cl:Phi}.
Then $\Psi(q,\xi_1,e_1)$ holds. Thus
$\exists x\exists \xi \exists e \Psi(x,\xi,e)$ also holds.
Since the parameters are all in $M$, we have
\begin{equation}
 M\models \exists x\ \exists \xi\  \exists e\  \Psi(x,\xi,e).
\end{equation} 
Thus there are $s,{\xi}_2, e_2\in M$
such that 
\begin{equation}\notag
\Phi(s,Z,q\restriction M)\land  (s\force \dot d_{{\xi}_2}=e_2)\land
e_2\subs a_s\land  (e_2\setm d)\subs a_s\setm a_{q\restriction M}.
\end{equation}
Since $e_2\setm d\subs a_s\setm a_{q\restriction M}$ and $(a_s\setm a_{q\restriction M})\cap Z\times Z=\empt$ by $\Phi(s,Z,q\restriction M)$, we have 
\begin{equation}
\notag 
(e_2\setm d)\cap Z\times Z=\empt
.
\end{equation}

Define the function 
\begin{equation}\notag
h:[a_s\setm a_{q\restriction M}, a_q\setm a_{q\restriction
    M}]\to \mc P({\pi}[a_s\cup a_q])
\end{equation}
by the formula
\begin{multline}
\label{eq:he} 
h\{x,y\}= {\pi}[a_s\cup a_q]\cap \min({\pi}(x), {\pi}(y))\cap\\
\bigcap\{X\in \mc A_q:x,y\in X\times X, \rank(X)\ge {\delta}\}.  
\end{multline}
So $h\{x,y\}$ is  as large as it is allowed by (\ref{eq:h}).

Then, by Lemma \ref{cl:Phi}, the condition
$r=\<a_r, f_r, \mc A_r\>$, where
$a_r=a_s\cup a_q$, $f_r=f_s\cup f_q\cup h$ and  
$\mc A_r=\mc A_s\cup \mc A_q$,
is a common extension of $q$ and $s$.

\begin{lemma}\label{lm:good}
$e_1$ and $e_2$ are good for $f_r$.  
\end{lemma}
\begin{proof}
We should check  
conditions (S1)-(S3).

Assume that $z\in e_1\cap e_2\cap ({\omega}_2\times {\omega})$,
 $x\in e_1\setm e_2\subs a_q\setm a_{q\restriction M}$, and  $y\in e_2\setm e_1\subs a_s\setm a_{q\restriction M}$.
Observe that $z,y\in M$, and so $\pi(z),\pi(y)\in M$ as well.

\smallskip
\noindent
{\bf (S1):}
Assume that  ${\pi}(z)< {\pi}(x), {\pi}(y)$.

We should show that 
$\pi(z)\in f_r\{x,y\}$.
However, $f_r\{x,y\}$ was defined by (\ref{eq:he}).
So we should show that
\begin{equation}\notag
\text{if $X\in \mc A_q$, $\rank(X)\ge {\delta}$, 
$x,y\in X\times X$ then 
${\pi}(z)\in X$.}  
\end{equation}
Since ${\pi}(y)\in M\cap  X\cap {\omega}_2$ 
and $\rank(M\cap {\lambda})\le \rank(X
)$
we have  $M\cap {\pi}(y)\subs X\cap {\pi}(y)$
by Lemma \ref{lm:main_lemma_7}(4).
Since ${\pi}(z)<{\pi}(y)$ and $\pi(z)\in M$ it follows that 
${\pi}(z)\in X$.

\smallskip
\noindent
{\bf (S2):}
Assume that  ${\pi}(z)<{\pi}(y)$.

We need to show that 
$f_r\{x,z\}\subs f_r\{x,y\}$.
Since $f_r\{x,z\}=f_q\{x,z\}$
and $f_r\{x,y\}$
was defined by (\ref{eq:he})
we should show that
\begin{equation}\notag
\text{if $X\in \mc A_q$, $\rank(X)\ge {\delta}$, 
$x,y\in X\times X$ then 
$f_q\{x,z\}\subs X$.}  
\end{equation}
Since ${\pi}(y)\in M\cap X$ we have 
${\pi}(z)\in M\cap {\pi}(y)\subs X\cap {\pi}(y)$ by Lemma \ref{lm:main_lemma_7}(4).

Since ${\pi}(z)\in X$,  $z\in {\omega}_2\times {\omega}$
and ${\omega}\subs X$ by (\ref{omega}),
it follows that $z\in X\times X$.
Since $x,z\in X\times X$ and $X\in \mc A_q$,
we have $f_q\{x,z\}\subs X$ by (\ref{eq:pdef}), which was to  be proved.

\smallskip
\noindent
{\bf (S3):}
Assume that  ${\pi}(z)<{\pi}(x)$.

We need to show that 
$f_r\{y,z\}\subs f_r\{x,y\}$.
Since $f_r\{y,z\}=f_s\{y,z\}$
and $f_r\{x,y\}$
was defined by (\ref{eq:he})
we should show that
\begin{equation}\notag
\text{if $X\in \mc A_q$, $\rank(X)\ge {\delta}$, 
$x,y\in X\times X$ then 
$f_s\{y,z\}\subs X$.}  
\end{equation}
Since
$z,y\in M$ we have
$f_s\{y,z\}\subs M$.

Moreover    $y\in X\times X$, and so
${\pi}(y)\in M\cap X\cap {\omega}_2$,
which implies
$M\cap \pi(y)\subs X\cap {\pi}(y)$ by Lemma \ref{lm:main_lemma_7}(4) .

Thus  $f_s\{y,z\}=f_s\{y,z\}\cap {\pi}(y)\subs M
\cap {\pi}(y)\subs X\cap {\pi}(y)\subs X$,
which was to be proved.
\end{proof}

Since $r\force \dot d_{\xi_1}=\check{e}_1\land \dot{d}_{\xi_2}=\check{e}_2
\land f\supset \check{f}_r$,
by Lemma \ref{lm:good} $r\force$ ``$ \dot d_{\xi_1}$ and 
$\dot{d}_{\xi_2}$ are good for $f$''. So  $f$
is a  $\Delta({\omega}_2\times {\lambda})$-function in $V[\mc G]$.

Since $|P|\le {\lambda}$ and so
${\left({\lambda}^{\oo}\right)}^{V[\mc G]}\le 
{\left((|P|+{\lambda})^{\oo}\right)}^{V}=
{\left({\lambda}^{\oo}\right)}^{V}=
{\lambda}$,
 the proof of Theorem \ref{tm:sDelta}
 is complete.
\end{proof}

\section{Space construction}

Assume that 
  $X$ is a scattered space. We say that a subspace   $Y\subs X$ is 
a {\em stem} of $X$ provided
\begin{enumerate}[(i)]
\item $ht(Y)=ht(X)$,
\item\label{stem2} $X\setm Y$ is closed discrete in	$X$.	
\end{enumerate}
Clearly  (\ref{stem2}) holds iff  
 every $x\in X$  has a neighborhood 
$U_x$ such that  
$U_x\setm \{x\}\subs Y$.

  \begin{proposition}\label{pr:stem}
  Assume that $X$ is an LCS space, $Y\subs X$ is a stem,
$SEQ(X)=\<{\kappa}_{\nu}:{\nu}<{\mu}\>$ and 
$SEQ(Y)=\<{\lambda}_{\nu}:{\nu}<{\mu}\>$.
Then
\begin{multline}
\{SEQ(Z): Y\subs Z\subs X\}=
\{ s\in {^{\mu}\text{Card}}: {\lambda}_\nu\le s({\nu})\le 
{\kappa}_\nu \text{ for each ${\nu}<{\mu}$} \}.
\end{multline}
  \end{proposition}

  \begin{proof}
Assume that $s\in {^{\mu}\text{Card}}$ such that ${\lambda}_{\nu}\le s({\nu})\le 
{\kappa}_\nu$ for each ${\nu}<{\mu}$.
For ${\nu}<{\mu}$  pick 
$Z_{\nu}\in [I_{\nu}(X)]^{s({\nu})}$ with $Z_\nu\supset I_\nu(Y)$.
Put $Z=\bigcup\{Z_{\nu}:{\nu}<{\mu}\}$.
Since $Y\subs Z$ and $Y$ is a stem, we have  
$I_{\nu}(Z)=  Z_{\nu}$ for ${\nu}<{\mu}$, and so $SEQ(Z)=s$.
  \end{proof}

\begin{theorem}\label{tm:space}
If there is a  $\Delta({\omega}_2\times{\lambda})$-function,
then there is a c.c.c poset $P$ such that in $V^P$
there is an $LCS$ space $X$ with stem $Y$ such that 
$SEQ(X)=\<{\lambda}\>_{{\omega}_2}$ and 
$SEQ(Y)=\<{\omega}\>_{{\omega}_2}$.  
\end{theorem}

\begin{corollary}\label{cor:space}
If there is a  $\Delta({\omega}_2\times{\lambda})$-function,
then there is a c.c.c poset $P$ such that in $V^P$   
every sequence $\mathbf s=\<s_{\alpha}:{\alpha}<{\omega}_2\>$ of infinite
cardinals with $s_{\alpha}\le {\lambda}$
is the cardinal sequence of some
locally compact scattered space.
\end{corollary}

\begin{proof}[Proof of Theorem \ref{tm:space}.]
Instead of constructing the topological space 
directly, we actually produce a certain ``graded poset'' which
guarantees the existence of the desired  locally compact scattered space.
We use the ideas from \cite{Bag} to formulate the properties of our required poset.

\begin{definition}\label{df:st-poset}
 Given two sequences $\mf t=\<{\kappa}_{\alpha}:{\alpha}<{\delta}\>$
and  $\mf s=\<\lambda_{\alpha}:{\alpha}<{\delta}\>$
 of
 infinite cardinals with $\lambda_\alpha\le \kappa_\alpha$, 
we say that a poset $\<T,\prec\>$ is a $\mf
 t$-{\em poset with an $\mf s$-stem}
iff the following conditions are satisfied:
\begin{enumerate}[(T1)]
\item $T=\bigcup\{T_{\alpha}:{\alpha}<{\delta}\}$ where
$T_{\alpha}=\{{\alpha}\}\times {\kappa}_{\alpha}$ for each $\alpha
<\delta$. 
Let $S_\alpha=\{\alpha\}\times \lambda_\alpha$, and 
$S=\bigcup\{S_\alpha:\alpha<\delta\}$.

\item For each $s\in T_{\alpha}$ and $t\in T_{\beta}$,
if $s \prec t$
  then ${\alpha}<{\beta}$ and $s\in S_\alpha$.

\item For every $\{s,t\}\in \br T;2;$ there is a finite subset
$i\{s,t\}$ of $S$ such that for each $u\in T$:

 \begin{equation}\notag
  (u\preceq s\land u\preceq t) \text{ iff }
 \text{$u\preceq v$ for some $v\in  i\{s,t\}$}.
  \end{equation}

\item For ${\alpha}<{\beta}<{\delta}$, if $t\in T_{\beta}$ then
the set $\{s\in S_{\alpha}:s\prec t\}$ is infinite.
\end{enumerate}
\end{definition}

\begin{lemma}\label{lm:Ba}
If there is a $\mf t$-poset with an $\mf s$-stem then there is an LCS space 
$X$ with stem $Y$ such that $SEQ(X)=\mf t$ and $SEQ(Y)=\mf s$.
\end{lemma}

Indeed, if  $\mc T=\<T,\prec\>$ is an $\mf s$-poset, we write
$U_{\mc T}(x)=\{y\in T:y\preceq x\}$ for $x\in T$, and  we denote
by $X_{\mc T}$ the topological space on $T$ whose subbase is the
family
\begin{equation}
 \{U_{\mc T}(x), T\setm U_{\mc T}(x):x\in T\},
\end{equation}
then $X_{\mc T}$ is our desired  LCS-space with stem. 

\vspace{1mm} So,  to prove Theorem \ref{tm:space} it will be
enough to show that a $\<\lambda\>_{\omega_2}$-poset
with an $\<\omega\>_{\omega_2}$-stem
may exist.

We follow the ideas of \cite{BS} to construct $P$.
Fix a  $\Delta({\omega}_2\times{\lambda})$-function  $f:\br  {\omega}_2\times{\lambda};2;\to \br \omega_2;<\omega; $.

\begin{definition}\label{df:bs}
Define the poset $\mc{P}=\<P,\preceq\>$ as follows.
The underlying set $P$ consists of  triples 
$p=\<a_p,\le_p,i_p\>$ satisfying the
following requirements:
\begin{enumerate}[(1)]
  \item $a_p\in \br {\omega}_2\times {\lambda};<{\omega};$,
\item $\le_p$ is a partial ordering on $a_p$ with the property that
if $x<_py$ then $x\in {\omega}_2\times {\omega}$ and ${\pi}(x)<{\pi}(y)$,
\item $i_p:\br a_p;2;\to \br a_p;<{\omega};$ is such that 
  \begin{enumerate}[(\arabic{enumi}.1)]
    \item if $\{x,y \}\in \br a_p;2;$ then 
      \begin{enumerate}[(\arabic{enumi}.\arabic{enumii}.1)]
          \item if $x,y\in {\omega}_2\times {\omega}$
and ${\pi}(x)={\pi}(y)$ then $i_p\{x,y\}=\empt$,
          \item if $x<_py$ then $i_p\{x,y\}=\{x\}$,
           \item if $x$ and $y$ are $<_p$-incomparable, then \\
               $i_p\{x,y\}\subs f\{x,y\}\times {\omega}$.
      \end{enumerate}
   \item if $\{x,y\}\in \br a_p;2;$ and  $z\in a_p$ then \\
         $\big((z\le_p x\land z\le_p y$)  { iff } 
         $\exists t\in i_p\{x,y\}\ z\le_p t\big).$
  \end{enumerate}
\end{enumerate} 
Set $p\preceq q$ iff $a_p\supseteq a_q$, $\le_p\restriction a_q=\le_q$ and
$i_p\restriction \br a_q;2;=i_q$.
\end{definition}

\begin{lemma}
$P$ satisfies ${\omega}_1$-c.c.. 
\end{lemma}

\begin{proof}
Let $\{p_{\nu}:{\nu}<{\omega}_1\}\subs P$, $p_v=\<a_{\nu}, \le_{\nu},
 i_{\nu}\>$. By thinning out our sequence we can assume that
 \begin{enumerate}[(i)]
 \item $\{a_{\nu}:{\nu}<{\omega}_1\}$ forms a $\Delta$-system
with kernel $a'$.
\item $i_{\nu}\restriction \br a';2;=i$.
\item $\le_{\nu}\restriction a'\times a'=\le$.
\item for each ${\nu}<{\mu}<\oo$ there is
a bijection ${\rho}_{{\nu},{\mu}}:a_{\nu}\to a_{\mu}$
such that 
\begin{enumerate}[(a)]
  \item ${\rho}_{{\nu},{\mu}}\restriction a'=\id$
\item $\pi(x)\le{\pi}(y)$ iff 
$ {\pi}({\rho}_{{\nu},{\mu}}(x))   \le {\pi}({\rho}_{{\nu},{\mu}}(y))$,
\item $x\le_{{\nu}} y$ iff ${\rho}_{{\nu},{\mu}}(x)\le_{\mu} {\rho}_{{\nu},{\mu}}(y)$,
\item  $x\in {\omega}_2\times {\omega}$ iff ${\rho}_{{\nu},{\mu}}(x)\in
  {\omega}_2\times {\omega}$,
\item ${\rho}_{{\nu},{\mu}}[i_{\nu}\{x,y\}]=i_{\mu}\{{\rho}_{{\nu},{\mu}}(x), 
{\rho}_{{\nu},{\mu}}(y)\}$.
\end{enumerate}
 \end{enumerate}
Now it follows from condition (3.1) and condition (iv) that if 
$\nu<\mu<\omega_2$ and $\{x,y\}\in \br a';2;$ then $i_\nu\{x,y\}=i_\mu\{x,y\}$.

Since $f$ is a $\Delta({\omega}_2\times {\lambda})$-function
there is ${\nu}<{\mu}<{\omega}_1$ such that 
$a_\nu$ and $a_\mu$ are good for $f$, i.e. (S1)--(S3)  hold. 
Define $r=\<a,\le, i \>$ as follows:
\begin{enumerate}[(a)]
\item $a=a_{\nu}\cup a_{\mu}$,
\item $x\le y$ iff $x\le_{\nu} y$ or $x\le_{\mu} y$
or there is $s\in a_{\nu}\cap a_{\mu}$ such that 
$x\le _{\nu} s\le_{\mu} y$ or $x\le_{\mu}s\le_{\nu} y$, 
\item $i\supset i_{\nu}\cup i_{\mu}$,
\item for $x\in a_{\nu}\setm a_{\mu}$ and 
$y\in a_{\mu}\setm a_{\nu}$, if $x$ and $y$ are $\le$-incomparable then 
  \begin{equation}\label{eq:i}
i\{x,y\}=\big(f\{x,y\}\times {\omega}\big)\cap\{t\in a: t\le x\land t\le
y\}.     
  \end{equation}
\item for $\{x,y\}\in \br a;2;$ with $x<y$, $i\{x,y\}=\{x\}$. 
\end{enumerate}
We claim that $r\in P$.
 
By the construction, we have $\le\restriction a_{\nu}\times
a_{\nu}=\le _{\nu}$ and $\le\restriction a_{\mu}\times
a_{\mu}=\le _{\mu}$.

\medskip
\noindent
{\bf Claim:} $\le$ is a partial order. 

We should check only the transitivity.
Assume $x\le y\le z$. If $x\le_{\nu} y\le_{\nu} z $ or 
$x\le _{\mu} y\le _{\mu}z$ then we are done.
Assume that $x\le_{\nu} u \le_{\mu} y \le_{\mu} z $ for some 
$u\in a_\nu\cap a_\mu$.
Then $x\le_{\nu} u  \le_{\mu} z $ so $x\le z$.

If $x\le_{\nu} u \le_{\mu} y  \le_{\mu} t\le_{\nu} z $
for some  $u,t\in a_{\nu}\cap a_{\mu}$, then   
$u\le _{\mu} t$,  which implies $u\le_{\nu} t$.
Thus $x\le_{\nu}u\le_{\nu} t \le_{\nu} z $ and so 
$x\le_{\nu} z$, and hence $x\le z$.

The other cases are similar to these ones.
\medskip 

(3.1.3) holds by the construction of $i$.

To show that $p$ is a condition we should finally check 
(3.2). Let $x,y\in a$ be $\le$-incomparable elements.
It is clear that if $u\le t$ for some $t\in  i\{x,y\}$ then $u\le x$ and $\le y$.
So we should check that 
\begin{enumerate}[(1)]
\item[$(*)$] if $z\le x$ and $z\le y$ then there is $t\in i\{x,y\}$
such that $z\le t$. 
\end{enumerate}

If $x,y,z\in a_{\nu}$ or $x,y,z\in a_{\mu}$ then it is clear
because $p_\nu,p_\mu \in P$.

\noindent
{\bf Case 1.} {\em $x,y\in a_{\nu}$ and 
$z\in a_{\mu}\setm a_{\nu}$.}
 
\noindent
{\bf Subcase 1.1} {\em $x,y\in a_{\nu}\setm a_{\mu}$.}

There are $x',y'\in a_{\nu}\cap a_{\mu}$ such that 
$z\le_{\mu} x'\le_{\nu} x$ and 
$z\le_{\mu} y'\le_{\nu} y$.
Then there is $t'\in i_{\mu}\{x',y'\}$ such that
$z\le_{\mu} t'$. Then $t'\in a_{\nu}\cap a_{\mu}$,
so $t'\le_{\nu} x,y$. Thus there is $t\in i_{\nu}\{x,y\}$ such that 
$t'\le_{\nu} t $, and so $z\le t$.
Since $i\{x,y\}=i_{\nu}\{x,y\}$, we are done.  

\noindent
{\bf Subcase 1.2} {\em $x\in a_{\nu}\setm a_{\mu}$ and 
$y\in a_{\nu}\cap a_{\mu}$}

Put $y'=y$, then 
proceed as in Subcase 1.1. 

\noindent
{\bf Case 2.} {\em $x,z\in a_{\nu} \setm a_{\mu}$ and 
$y\in a_{\mu}\setm a_{\nu}$.}

Then $z\le_{\nu} y'\le_{\mu} y$ for some $y'\in a_{\nu}\cap a_{\mu}$.
Then there is $t\in i_{\nu}\{x,y'\}$ such that 
$z\le_{\nu} t$. Clearly $t\le x,y$.  We show that $t\in i\{x,y\}$.

If $t=y'$ then $t\le x,y$ and $\pi(t)\in f\{x,y\}$ by
(S1). Thus $t\in i\{x,y\}$.  

Assume that $t<_{\nu} y'$.
Then $\pi(t)\in f\{x,y'\}\subs f\{x,y\}$ by (S2), because
$y'\in a_{\nu}\cap a_{\mu}$ and ${\pi}(y')<{\pi}(y)$. 
Thus $t\in i\{x,y\}$ by (\ref{eq:i}).
\end{proof}

Assume that $\mc G$ is a $\mc P$-generic filter. We claim that  
if we take
\begin{equation}\notag
\ll=\cup\{\le_p:p\in \mc G\}. 
\end{equation} 
then $\<\omega_2\times \lambda,\ll\>$ is a 
$\<\lambda\>_{\omega_2}$-poset 
with an $\<\omega\>_{\omega_2}$-stem.
By standard density arguments, $\ll$ is a partial order on $\omega_2\times \lambda$
which satisfies (T4). Moreover, every $p\in P$ satisfies (2), so (T2) also holds.
Finally the function 
\begin{equation}\notag
\mathop i=\bigcup\{i_p:p\in \mc G\} 
\end{equation} 
witnesses (T3) because every $p\in P$ satisfies (3.2).
\end{proof}

\bigskip

\noindent{\bf Acknowledgement}
I would like to thank Professor Juan Carlos Martinez  for several useful discussions, comments and suggestions.

\end{document}